\documentclass[11pt]{article}

\usepackage{latexsym,mathrsfs}
\usepackage{amsmath,amssymb}
\usepackage{amsthm,enumerate,verbatim}
\usepackage{amsfonts}
\usepackage{graphicx}
\usepackage{algorithm}
\usepackage{algorithmic}
\usepackage{url}
\usepackage{color} 
\usepackage{appendix}
\usepackage{lscape}

\definecolor{brightpink}{rgb}{1.0, 0.0, 0.5}

\newtheorem{theorem}{Theorem}

\newtheorem{definition}{Definition}
\newtheorem*{definition*}{Definition}

\newtheorem{remark}{Remark}

\numberwithin{equation}{section}
\numberwithin{table}{section}
\numberwithin{figure}{section}
\def \R{{\mathbb R}}
\def \C{{\mathbb C}}
\def \H{{\mathbb H}} 
\def\real{\mathop{\mathrm{Re}}}
\def\imag{\mathop{\mathrm{Im}}}

\newcommand {\mat}  [1] {\left[\begin{array}{#1}}
\newcommand {\rix}      {\end{array}\right]}

\def\real{\mathop{\mathrm{Re}}}
\def\imag{\mathop{\mathrm{Im}}}
\newcommand{\eproof}{\space
    {\ \vbox{\hrule\hbox{\vrule height1.3ex\hskip0.8ex\vrule}\hrule}}\par}

\title{Characterizing matrices with eigenvalues in an LMI region: \\ A dissipative-Hamiltonian approach}
\date{}

\author{Neelam Choudhary\thanks{
 School of Engineering and Applied Sciences, Department of Mathematics, Bennett University, Greater Noida-201310, Uttar Pradesh, India; \texttt{neelam.choudhary@bennett.edu.in. }}
 \qquad Nicolas Gillis\thanks{
 Department of Mathematics and Operational Research,
Facult\'e Polytechnique, Universit\'e de Mons, Rue de Houdain~9, 7000 Mons, Belgium; \texttt{nicolas.gillis@umons.ac.be}. N.\ Gillis  acknowledges the support by the Fonds de la Recherche Scientifique - FNRS and the Fonds Wetenschappelijk Onderzoek - Vlanderen (FWO) under EOS Project no O005318F-RG47, and by the Francqui Foundation. 
}
 \qquad Punit Sharma\thanks{Department of Mathematics, Indian Institute of Technology Delhi, Hauz Khas, New Delhi-110016, India; \texttt{punit.sharma@maths.iitd.ac.in}.
P.\ Sharma acknowledges the support of the DST-Inspire Faculty Award (MI01807-G) by the Government of India and Institute SEED Grant (NPN5R) by IIT Delhi.
 }
}

\begin{document}

\maketitle

\begin{abstract}
In this paper, we provide a dissipative Hamiltonian (DH) characterization for the set of matrices whose eigenvalues belong to a given LMI region. 
This characterization is a generalization of that of Choudhary et al.\ 
(Numer.\ Linear Algebra Appl, 2020) to any LMI region. 
It can be used in various contexts, which we illustrate on the nearest $\Omega$-stable matrix problem: given an LMI region $\Omega \subseteq \C$ and a matrix $A \in \C^{n,n}$, find the nearest matrix to $A$ whose eigenvalues belong to $\Omega$. 
Finally, we generalize our characterization to more general regions that can be expressed using LMIs involving complex matrices. 
\end{abstract}


\section{Introduction}

In this paper, we study matrices $A\in \R^{n,n}$ whose eigenvalues belong to a subset of the complex plane, $\Omega \subseteq \C$, such matrices are called $\Omega$-stable. 
\begin{definition}{\rm ($\Omega$-stability)}
For $\Omega \subseteq \mathbb C$, the matrix $A \in \mathbb R^{n,n}$ is said to 
 be $\Omega$-stable if every eigenvalue of $A$ lie inside the region $\Omega$. 
\end{definition}
The two most famous examples of $\Omega$-stable matrices are Hurwitz stable matrices for which  {$\Omega=\{z\in \C:~\real{z} < 0\}$}, and Schur stable matrices for which $\Omega=\{z\in \C:~|z| < 1\}$. 
Hurwitz stable matrices play a significant role in the study of 
linear time-invariant (LTI) systems of the form
\begin{eqnarray*}
\dot{x}(t)=A x(t)+Bu(t), \quad y(t)=Cx(t),
\end{eqnarray*}
where, for all $t\in \R$, $x(t) \in \R^n$, $u(t) \in \R^{m}$, 
$y(t)\in \R^{p}$, $A\in \R^{n,n}$, $B\in \R^{n,m}$, and $C\in \R^{p,n}$. In fact, such a system is stable if $A$ is Hurwitz stable. 
Moreover, the transient response of such a system is  directly related to the location of its poles~\cite{ChilG96} in the complex plane. The poles in a specific region in the complex plane 
can bound the maximum overshoot, the frequency of oscillatory modes, the delay time, the rise time, and the settling time. The problem of locating all the closed-loop poles of a controlled system inside a specific region $\Omega \subseteq \C$ is known as the \emph{$\Omega$-pole placement problem} and has appeared in several applications~\cite{ChilG96,Had92,Bha09,WanF06,ChiGA99,YaoRH14}. 

For these reasons, characterizing $\Omega$-stable matrices is an important topic in numerical linear algebra and control. 
In this paper, we focus on regions of the complex plane that can be expressed by linear matrix inequalities (LMIs). 
Such sets are referred to as \emph{LMI regions}~\cite{ChilG96} and defined as follows. 
\begin{definition}[LMI Region, \cite{ChilG96}] \label{def:lmireg} 
A subset $\Omega \subseteq \C$ is called an LMI region if there exists a symmetric matrix $B \in \R^{s,s}$ 
and a matrix $C \in \R^{s,s}$ such that 
\begin{equation}\label{eq:eqlmidef1}
\Omega =\left \{z\in \C:~f_\Omega(z) \prec 0 \right \},
\end{equation}
where 
\begin{equation}\label{eq:eqlmidef2}
f_\Omega 
\colon 
\mathbb{C} \mapsto \H^{s,s} \; \text{ is given by } \; z \mapsto f_\Omega(z) := B + z C +\overline z C^T ,  
\end{equation}
where 
$\H^{s,s}$ is the set of Hermitian matrices with real eigenvalues, 
that is, 
$\H^{s,s}  = \{ X \in \C^{s,s} \ : \ X = X^* \}$ with $X^*$ the conjugate transpose of $X$, 
and for $X \in \H^{s,s}$, $X \prec 0$ means that $X$ is negative definite, that is, its eigenvalues are negative. 
\end{definition} 
\noindent The function $f_\Omega(z)$ is called the characteristic function of the LMI region $\Omega$, 
and $s$ is called the order of $f_\Omega(z)$. 
The characteristic function of an LMI region is 
not unique~\cite{Kus19}. 
Since $f_{\Omega}(\overline z)=(f_{\Omega}(z))^T$, any LMI region is symmetric along the real axis. 
An LMI region is convex, and so is the intersection of two or more LMI regions.  Due to the strict inequality ``$\prec$" in \eqref{eq:eqlmidef1}, the LMI regions are open. Furthermore, 
the LMI regions are invariant under congruence transformations. We refer to~\cite{Kus19} for other topological and geometrical properties of the LMI regions. 
A large number of subsets in the complex plane can be expressed as LMI regions; for example, conic sectors, vertical half-planes, vertical strips, discs, horizontal strips, ellipses, parabolas, and  hyperbolic sectors; see~\cite{BisPT21} and 
Section~\ref{sec:specificlmis}. 
The set of LMI regions is dense in the set of convex regions symmetric to the real axis, which are relevant for control systems~\cite{ChilG96,BisPT21}.

In~\cite{GilS17}, authors characterized $\Omega$-stable matrices using the so-called dissipative Hamiltonian (DH) matrices. 
\begin{definition}[DH matrix]
A matrix $A \in \R^{n,n}$ is said to be a DH matrix if $A=(J-R)Q$ for some 
$J,R,Q \in \R^{n,n}$ such that $J^T =-J$, $R\succeq 0$, and $Q \succ 0$.
\end{definition} 
A DH matrix is always Hurwitz stable, that is, all its  eigenvalues are in the left half of the complex plane~\cite{GilS17}. 
The term DH is inspired by the DH systems in which the state matrix has the form $A = (J-R)Q$, where $J^T=-J$  is the structure matrix describing the flux among energy storage elements, $R$ is a positive semidefinite matrix describing the energy dissipation in the system, and $Q$ is a positive definite matrix that describes the energy of the system~\cite{Sch06,SchM95}.  
By replacing the constraint on $R$, namely $R\succeq 0$, by other LMI constraints on $(J,R,Q)$, DH matrices can represent different types of $\Omega$-stable matrices. 
This was studied in~\cite{ChoGS20} where $\Omega$-stable matrices were written as DH matrices where $\Omega$ could be vertical strips, disks, conic sectors, and their intersection; see the next section for more details. 
An application of these characterizations is to solve the nearest $\Omega$-stable matrix problem. 
For example, in system identification, one needs to identify a $\Omega$-stable system from observations~\cite{OrbNV13}. In fact, sometimes numerical or modelling errors or approximation processes may 
produce an unstable system in place of a stable one. The unstable system then has to be approximated by a nearby stable system without perturbing its entries too much. More precisely, for a region $\Omega \in \C$ and a matrix $A\in \R^{n,n}$, this requires to solve the following optimization problem 
\begin{equation}\label{eq:probdef}
\inf_{X \in \mathcal S_{\Omega}} {\|A-X\|}_F^2, 
\end{equation} 
where $\|\cdot\|_F$ stands for the Frobenius norm and $\mathcal S_{\Omega}$ is the set of all $\Omega$-stable matrices. 
The DH characterization of stable matrices has been proven very effective in solving several nearness problems for LTI control systems.
For example, distance to $\Omega$-stability~\cite{GilS17,ChoGS20, mamakoukas2020learning}, nearest admissible descriptor system problem~\cite{GilMS17}, distance to passivity~\cite{GilS17b},  minimal-norm-static-state feedback problem~\cite{GilS21}, and 
learning data-driven stable Koopman operators~\cite{mamakoukas2020learning2}.  
This DH characterization was also used recently to design an optimization-based algorithm for parametric model order reduction of LTI dynamical systems~\cite{schwerdtner2022structured}.  


\paragraph{Contribution and outline of the paper} 

This paper is organized as follows. 
In Section~\ref{sec:DHcarac}, for a given LMI region $\Omega \subseteq \C$, we characterize the set of all $\Omega$-stable matrices as DH matrices of the form $A = (J-R)Q$ with LMI constraints on the triplet $(J,R,Q) \in (\R^{n,n})^3$. This characterization generalizes the work on~\cite{ChoGS20} that only considered three types of LMI regions, namely conic sectors, vertical strips and disks.  In Section~\ref{sec:specificlmis}, we provide several examples of such LMI regions, including parabolas, ellipsoids, hyperbolas and horizontal strips; this is the first time DH characterizations of such regions are given via LMIs on the triplet $(J,R,Q)$. 
In Section~\ref{sec:appl}, we illustrate the use of these characterizations to solve the nearest $\Omega$-stable matrix problem~\eqref{eq:probdef}. 
Finally, in Section~\ref{sec:extebdedLMIs}, we extend our characterizations of LMI regions for complex matrices, where the set $\Omega$ is not necessarily symmetric with respect to the real line.

\paragraph{Notation} Throughout the paper, $X^T$ and
$\|X\|$  stand for the transpose and the spectral norm of a real matrix $X$, respectively.
We write $X\succ 0$ $(X \prec0)$ and $X\succeq 0$ $(X \preceq 0)$ if $X$ is symmetric and positive definite (negative definite)
or positive semidefinite (negative semidefinite), respectively.
By $I_m$ we denote the identity matrix of size $m \times m$. The  Kronecker  product is represented by $\otimes$ and we refer to~\cite{HorJ85} for the standard properties of the Kronecker product. The set of $n \times n$ Hermitian matrices is denoted by $\H^{n,n}$.

\section{DH characterization of matrices with eigenvalues in generic LMI regions} \label{sec:DHcarac}

In this section, we consider matrices with eigenvalues in some generic LMI regions, and for them, we provide a parametrization using the DH matrices. This will allow us for example in Section~\ref{sec:appl} to use standard optimization tools to find a nearby matrix to a given matrix with eigenvalues all in the given LMI region. 

The following result from~\cite{ChilG96} is crucial for our DH formulation of the set of 
$\Omega$-stable matrices, where $\Omega$ is an LMI region.  
\begin{theorem}{\rm \cite[Theorem 2.2]{ChilG96}} \label{Thm:lmicharmain}
 Let $\Omega$ be an LMI region given by~\eqref{eq:eqlmidef1} and let $A \in \R^{n,n}$. Then $A$ is $\Omega$-stable if and only if there exists a symmetric matrix $X \in \R^{n,n}$ such that $X \succ 0$ and 
\[
\mathcal M_{\Omega}(A,X):= B  \otimes X + C \otimes (AX) + 
C^T \otimes (AX)^T \prec 0.
\] 
\end{theorem}

Let us illustrate Theorem~\ref{Thm:lmicharmain} on the two most well-known examples.  
\begin{enumerate}
    \item For Hurwitz stability with $\Omega=\{z\in \C:~\real{z} <  0\}$, we simply take $B = 0$ and $C = 1$ to obtain  
    $f_\Omega(z) = \real{z}$. 
    This choice leads to the classical Lyapunov stability criterion, namely there exists $X \succ 0$ such that $AX + XA^T \prec 0$. 
    
    \item For Schur stability with $\Omega=\{z\in \C:~|z| < 1\}$,  we take $B = -\binom{1 \; 0}{0 \; 1}$ and 
    $C = \binom{0 \; 1}{0 \; 0}$, so that 
    \[
 f_\Omega(z) = 
 \mat{cc} 
 -1 & z \\ 
 \overline{z} & -1 
 \rix 
 = 
  \underbrace{\mat{cc} 
 -1 & 0 \\ 
  0 & -1 
 \rix}_{=B}
 + 
 z \underbrace{\mat{cc} 
 0 & 1 \\ 
 0 & 0
 \rix}_{=C} 
 + 
 \bar z \underbrace{\mat{cc} 
 0 & 0  \\ 
 1 & 0 
 \rix}_{=C^\top} . 
 \]
 In fact, $f_\Omega(z) \prec 0$ if and only if the trace of $f_\Omega(z)$ is negative (the sum of the eigenvalues is negative) and the determinant is positive (the product of the eigenvalues is positive). 
 The trace is equal to $-2$, 
 while the determinant is given by $1 - |z|^2$, which leads to the region $|z|^2 < 1$, as desired.  This choice leads to the well-known condition for discrete-time stability: 
 there exists $X \succ 0$ such that $\mat{cc} 
-X & AX \\ 
 X A^T & -X  
 \rix \prec 0$ which is equivalent, using the Schur complement, to 
 $X - AXA^T \succ 0$. 
\end{enumerate}

To characterize matrices with eigenvalues in any LMI region, 
we rephrase the definition of a DH matrix by relaxing the semidefinite constraint on $R$ and obtaining the following characterization of $\Omega$-stable matrices. 
\begin{theorem} \label{th:genJRQ}
Let $ \Omega \subseteq C$ be an LMI region defined by~\eqref{eq:eqlmidef1} and let $B,C \in \R^{n,n}$ be the corresponding matrices as given in~\eqref{eq:eqlmidef2}. 
Consider a matrix $A \in \R^{n,n}$. Then $A$ is $\Omega$-stable if and only if $A=(J-R)Q$ for some 
$J,R,Q \in \R^{n,n}$ such that $J^T =-J$, $R^T=R $, $Q \succ 0$, 
and 
\begin{equation} \label{eq:thm2}
\mathcal M_{\Omega}(J,R,Q) 
: = B \otimes Q^{-1} +  (C-C^T ) \otimes J - (C + C^T ) \otimes R 
\quad \prec \quad 0 . 
\end{equation}
\end{theorem}

\proof First suppose that $A=(J-R)Q$ for some 
$J,R,Q \in \R^{n,n}$ such that $J^T =-J$, $R^T=R $, $Q \succ 0$, and $\mathcal M_{\Omega}(J,R,Q) \prec 0$. Let $\lambda \in \C$
be an eigenvalue of $A$ and $v \in \C^n \setminus \{0\}$ be the corresponding eigenvector, that is, 
\begin{equation*}
Av =(J-R)Qv=\lambda v.
\end{equation*}
Thus
\begin{eqnarray}
\left (B+ \lambda C + \overline \lambda C^T  \right)v^*Qv  &=& B \otimes v^*Qv + \lambda C \otimes v^*Qv + \overline \lambda C^T  \otimes v^*Qv  \nonumber \\ 
&=&  (I\otimes v)^* \left ( B \otimes Q^{-1} +  (C-C^T ) \otimes J - (C + C^T ) \otimes R\right )  (I \otimes v) \nonumber  \\
&=&  (I\otimes v)^* \mathcal M_{\Omega}(J,R,Q) (I \otimes v) \nonumber \\
&\prec & 0, \nonumber 
\end{eqnarray} 
since $\mathcal M_{\Omega}(J,R,Q) \prec 0$ and $Q \succ 0$.  
This implies that $(B+ \lambda C + \overline \lambda C^T ) \prec 0$, and hence from~\eqref{eq:eqlmidef1} $\lambda \in \Omega$.
Conversely, let $A$ be $\Omega$-stable. Then by Theorem~\ref{Thm:lmicharmain}, there exists $X \succ 0$ such that 
\begin{equation}\label{thm:eqgen1}
\mathcal M_{\Omega}(A,X)=B  \otimes X + C \otimes (AX) + C^T  \otimes (AX)^T  \prec 0.
\end{equation}
By setting 
\[
Q=X^{-1},\quad R=-\frac{(AX)+(AX)^T }{2}, \quad \text{and}\quad J=\frac{(AX)-(AX)^T }{2},
\]
we have $J^T =-J$, $R ^T=R$, and $Q \succ 0$ so that $A=(J-R)Q$. Also, in view of~\eqref{thm:eqgen1}, we have that 
\begin{eqnarray*}
\mathcal M_{\Omega}(J,R,Q) &=&   B \otimes Q^{-1} +  (C-C^T ) \otimes J - (C + C^T ) \otimes R \\
&=& B \otimes X +  (C-C^T ) \otimes \frac{(AX)-(AX)^T }{2} + (C + C^T ) \otimes \frac{(AX)+(AX)^T }{2} \\
&=& B  \otimes X + C \otimes (AX) + C^T  \otimes (AX)^T  \\
&\prec & 0.
\end{eqnarray*}
This completes the proof.
\eproof

If we apply Theorem~\ref{th:genJRQ} to Hurwitz and Schur stable matrices, we obtain the following DH characterizations: 
\begin{enumerate}
    \item For Hurwitz stability, plugging $B = 0$ and $C = 1$ in \eqref{eq:thm2} implies that $A = (J-R)Q$ such that $J^T =-J$, $R^T=R $, $Q \succ 0$ is Hurwitz stable if and only if $R \succ 0$, as shown in~\cite{GilS17}.  
    
    \item For Schur stability, plugging $B = -\binom{1 \; 0}{0 \; 1}$ and 
    $C = \binom{0 \; 1}{0 \; 0}$ in \eqref{eq:thm2} implies that $A = (J-R)Q$ such that $J^T =-J$, $R^T=R $, $Q \succ 0$ is Schur stable if and only if 
    \[ 
 \mat{cc} Q^{-1} & -J+R 
 \\ J+R & Q^{-1}
\rix \succ 0, 
 \]
  as shown in~\cite{ChoGS20}. 
\end{enumerate}

In the next section, we apply Theorem~\ref{th:genJRQ} in the same way for other LMI regions.

\section{Special LMI regions and DH parametrization} \label{sec:specificlmis} 

We have provided in Theorem~\ref{th:genJRQ} the DH parametrization of matrices with eigenvalues inside a generic LMI region. As a corollary of which, in this section, we provide LMI constraints on the matrix triplets $(J,R,Q)$ of a DH matrix for some particular regions. In fact, as soon as an LMI region is defined via the matrices $B$ and $C$, see Definition~\ref{def:lmireg}, we can provide its corresponding DH characterization in terms of DH matrices. 
In this section, we provide the following examples: 
 left and right conic sectors, 
 disks centered on the real line, 
 vertical left and right halfplanes, 
 ellipsoid centred on the real line, 
 left and right parabolic regions centred on the real line, 
 left and right hyperbolas with vertices on the real line, 
 and horizontal strip. By following~\cite[Table III]{BisPT21}, we first define these regions as LMI regions for $s=2$, and, in Table~\ref{tab:specialdhlmi}, we provide the DH characterization of matrices with eigenvalues in these regions. 
 Here are the considered LMI regions: 


 
\begin{itemize}

\item Left Conic sector: the left conic sector region of parameters $a,\,\theta \in \R$ with $0 \leq \theta \leq \pi/2$, denoted by $\Omega_{C_L}(a,\theta)$, is defined as
\[
\Omega_{C_L}(a,\theta)
:=
\left\{ x+iy \in \C \ \big| \ \sin(\theta) (x-a) < \cos(\theta) y < -\sin(\theta) (x-a),~ x \leq a \right\}.
\]
The conic sector $\Omega_{C_L}(a,\theta)$ can be characterized in form of~\eqref{eq:eqlmidef1} of an LMI region with matrices
\[
B=\mat{cc} -a \sin(\theta) & 0 \\ 0 & -a \sin(\theta) \rix \quad 
\text{and}\quad C=\frac{1}{2}\mat{cc} \sin(\theta) &\cos(\theta)\\-\cos(\theta) & \sin(\theta) \rix.
\]

\item Right Conic sector: the right conic sector region of parameters $a,\,\theta \in \R$ with $0 \leq \theta \leq \pi/2$, denoted by $\Omega_{C_R}(a,\theta)$, is defined as
\[
\Omega_{C_R}(a,\theta)
:=
\left\{ x+iy \in \C \ \big| \ \sin(\theta) (a-x) < \cos(\theta) y < -\sin(\theta) (a-x),~ x > a \right\}.
\]

The conic sector $\Omega_{C_R}(a,\theta)$ is an LMI region in form of~\eqref{eq:eqlmidef1} with matrices
\[
B=\mat{cc} a\sin(\theta)  & 0 \\ 0 & a \sin(\theta) \rix \quad 
\text{and}\quad C=\frac{1}{2}\mat{cc} -\sin(\theta) &\cos(\theta)\\-\cos(\theta) & -\sin(\theta) \rix.
\]

\item Disks centred on the real line:
the disk centred  at $(q,0)$ with radius $r>0$, denoted by $\Omega_D(q,r)$, is defined as
\[
\Omega_D(q,r)
:= \left\{ z\in \C \ \big| \ |z-q|<r\right\}.
\]
The disk region $\Omega_D(q,r)$ can be characterized in form of~\eqref{eq:eqlmidef1} of an LMI region with matrices
\[
B=\mat{cc} -r & q \\q &-r \rix \quad \text{and}\quad 
C=\mat{cc}0&0\\-1&0\rix.
\]

\item Vertical strip: the vertical strip region of parameters $ h < k$, denoted by $\Omega_V(h,k)$, is defined as
\[
\Omega_V(h,k)
:= \left\{ x+iy \in \C \ \big| \ h < x < k \right\}.
\]
The vertical strip region $\Omega_V(h,k)$ can be characterized in form of~\eqref{eq:eqlmidef1} of an LMI region with matrices
\[
B=\mat{cc} -k & 0\\0 & h \rix \quad \text{and} \quad
C=\frac{1}{2} \mat{cc} 1 &0 \\0 & -1 \rix.
\]

Note that $h$ (resp.\@ $k$) can possibly be equal to $-\infty$ (resp.\@ $+\infty$)  in which case $\Omega_V$ is a half space.
In particular, $\Omega_V(0,+\infty)$ is the open left half of the complex plane, corresponding to stable matrices for continuous LTI systems.

More precisely, the vertical halfplane left of $k$, that is, 
$\Omega_V(-\infty,k)
:= \left\{ x+iy \in \C \ \big| \ x < k \right\}$ is an LMI region with matrices $
B=\mat{cc} -k & 0\\0 & -1 \rix$ and $ C=\frac{1}{2} \mat{cc} 1 &0 \\0 & 0 \rix $. 

Similarly, the vertical halfplane right of $h$, that is, 
$\Omega_V(h,\infty)
:= \left\{ x+iy \in \C \ \big| \ h< x  \right\}$ is an LMI region with matrices $
B=\mat{cc} h & 0\\0 & -1 \rix$ and $ C=\frac{1}{2} \mat{cc} -1 &0 \\0 & 0 \rix $.

\item Ellipsoid:~ the ellipse centred at $(q_e,0)$ with horizontal radius $a_e$ and vertical radius $b_e$,  denoted by $\Omega_E(q_e,a_e,b_e)$, is defined as
\[
\Omega_E(q_e,a_e,b_e)
:=
\left\{ x+iy \in \C \ \big| \ \frac{(x-q_e)^2}{a_{e}^2}+\frac{y^2}{b_{e}^2} < 1 \right\}.
\]
The ellipsoid $\Omega_E(q_e,a_e,b_e)$ can be characterize in form of~\eqref{eq:eqlmidef1} of an LMI region with matrices $B= \mat{cc}-2a_h & -2q_e \\ -2q_e & -2a_e \rix$ and $C=\mat{cc} 0 & (1+\frac{a_e}{b_e})\\ (1-\frac{a_e}{b_e}) & 0 \rix $.

\item Left parabolic region: the parabolic region symmetric to $x$-axis with centred at $(q_p,0)$, and curvature $c_p > 0$, denoted by $\Omega_{P_L}(q_p,c_p)$, defined as 
\[
\Omega_{P_L}(q_p,c_p)
:=
\left\{ x+iy \in \C \ \big| \ y^2 < -\frac{2}{c_p}(x-q_p)  \right\}.
\]
The parabolic region $\Omega_{P_L}(q_p,c_p)$ can be characterized in form of~\eqref{eq:eqlmidef1} of an LMI region with matrices
\[
B=\mat{cc} -1 & 0 \\ 0 & -q_p \rix \quad \text{and} \quad
C=\frac{1}{2}\mat{cc} 0 & \sqrt{\frac{c_p}{2}} \\ -\sqrt{\frac{c_p}{2}} & 1
\rix.
\]

\item Right parabolic region: the parabolic region symmetric to $x$-axis with centred at $(q_p,0)$, and curvature $c_p > 0$, denoted by $\Omega_{P_R}(q_p,c_p)$, defined as 
\[
\Omega_{P_R}(q_p,c_p)
:=
\left\{ x+iy \in \C \ \big| \ y^2 < \frac{2}{c_p}(x-q_p)  \right\}.
\]
The parabolic region $\Omega_{P_R}(q_p,c_p)$ can be characterized in form of~\eqref{eq:eqlmidef1} of an LMI region with matrices
\[
B=\mat{cc} -1 & 0 \\ 0 & q_p \rix \quad \text{and} \quad
C=\frac{1}{2}\mat{cc} 0 & \sqrt{\frac{c_p}{2}} \\ -\sqrt{\frac{c_p}{2}} & -1
\rix.
\]

\item Left hyperbolic region: the left hyperbolic region with the semi-major axis $a_h$, $a_h > 0$ and semi-minor axis $b_h$, $b_h > 0$ and two vertices at $(a_h,0)$ and $(-a_h,0)$, is denoted by $\Omega_{Hyp,L}(a_h,b_h)$, defined as 
\[
\Omega_{Hyp,L}(a_h,b_h)
:=
\left\{ x+iy \in \C \ : \ x<0,~\frac{x^2}{a_h^2}-\frac{y^2}{b_h^2} - 1 > 0 \right\}.
\]
The hyperbolic region $\Omega_{Hyp,L}(a_h,b_h)$ can be characterized in form of~\eqref{eq:eqlmidef1} of an LMI region with matrices
\[
B= \mat{cc} 0 & 1 \\ 1 & 0 \rix \quad \text{and}\quad
C= \mat{cc}\frac{1}{2a_h} & \frac{1}{2b_h} \\ -\frac{1}{2b_h} & \frac{1}{2a_h} \rix.
\]

\item Right hyperbolic region: the right hyperbolic region with the semi-major axis $a_h$ $a_h > 0$ and semi-minor axis $b_h$, $b_h > 0$ and two vertices at $(a_h,0)$ and $(-a_h,0)$ is denoted by $\Omega_{Hyp,R}(a_h,b_h)$, defined as
\[
\Omega_{Hyp,R}(a_h,b_h)
:=
\left\{ x+iy \in \C \ : \ x>0,~\frac{x^2}{a_h^2}-\frac{y^2}{b_h^2} - 1 > 0 \right\}.
\]
The hyperbolic region $\Omega_{Hyp,R}(a_h,b_h)$ can be characterized in form of~\eqref{eq:eqlmidef1} of an LMI region with matrices
\[
B= \mat{cc} 0 & 1 \\ 1 & 0 \rix \quad \text{and}\quad
C= \mat{cc}-\frac{1}{2a_h} & \frac{1}{2b_h} \\ -\frac{1}{2b_h} & -\frac{1}{2a_h} \rix.
\]

\item Horizontal strip: the horizontal strip symmetric to the real axis with hight $w$, denoted by $\Omega_H(w)$, defined as
\[
\Omega_H(w)
:=
\left\{ x+iy \in \C \ : |y| < w\right\}.
\]
The horizontal strip $\Omega_{Hyp,L}(a_h,b_h)$ can be characterized in form of~\eqref{eq:eqlmidef1} of an LMI region with matrices
\[
B= \mat{cc}  -w & 0 \\ 0 & -w \rix \quad \text{and} \quad 
C= \frac{1}{2} \mat{cc} 0 & 1 \\ -1 & 0 \rix.
\]
\end{itemize}

We note that DH characterization for the sets of $\Omega$-stable matrices was obtained in~\cite{ChoGS20}, where $\Omega \in \{\Omega_{C_L}(a,\theta), \Omega_D(q,r), \Omega_V(h,k)\}$. 
Because of Theorem~\ref{th:genJRQ}, we can now extend these results for other LMI regions as well, see Table~\ref{tab:specialdhlmi}. 


\begin{table}[ht]
    \centering
\begin{tabular}{|c|c|}
\hline
\textbf{LMI Region} & \textbf{Constraints on $J^T=-J$, $R^T=R$, and $Q\succ 0$} \\ \hline
Left Conic Sector: $\Omega_{C_L}(a,\theta)$ 
& $\mat{cc} \sin(\theta) (aQ^{-1}+R) & -\cos(\theta) J \\ \cos(\theta) J & \sin(\theta) (aQ^{-1}+R) \rix \succ 0 $\\ \hline
Right Conic Sector: $ \Omega_{C_R}(a,\theta)$ & $\mat{cc} -\sin(\theta) (aQ^{-1}+R) & -\cos(\theta) J \\ \cos(\theta) J & -\sin(\theta) (aQ^{-1}+R) \rix \succ 0 $ \\ \hline
Disk: $\Omega_D(q,r)$ & $\mat{cc} rQ^{-1} & qQ^{-1}-J+R \\ qQ^{-1}+J+R & rQ^{-1}
\rix \succ 0 $ \\ \hline
Vertical strip: $ \Omega_V(h,k)$ & $\mat{cc} kQ^{-1}+R & 0 \\ 0 & -hQ^{-1}-R
\rix \succ 0$\\ \hline
Left halfplane: $ \Omega_V(-\infty,k)$ & $ kQ^{-1}+R \succ 0$ \\ \hline
Right halfplane: $ \Omega_V(h,\infty)$ & $ -hQ^{-1}-R \succ 0$\\ \hline
Ellipsoid: $ \Omega_E(q_e,a_e,b_e)$ & $\mat{cc} a_eQ^{-1} & q_eQ^{-1} -  \frac{a_e}{b_e} J + R \\ q_eQ^{-1} +  \frac{a_e}{b_e} J + R & a_eQ^{-1}\rix \succ 0 $ \\ \hline
Left parabolic region: $ \Omega_{P_L}(q_p,c_p)$ & $\mat{cc} Q^{-1} & -\sqrt{\frac{c_p}{2}}J \\ \sqrt{\frac{c_p}{2}}J  & q_pQ^{-1}+ R \rix \succ 0$ \\ \hline
Right parabolic region: $ \Omega_{P_R}(q_p,c_p)$ & $\mat{cc} Q^{-1} & -\sqrt{\frac{c_p}{2}}J \\ \sqrt{\frac{c_p}{2}}J  & -q_pQ^{-1}- R \rix \succ 0$\\ \hline
Left hyperbola: $ \Omega_{Hyp,L}(a_h,b_h)$ & $ \mat{cc} \frac{R}{a_h} & -Q^{-1} - \frac{J}{b_h}\\ -Q^{-1}+\frac{J}{b_h} & \frac{R}{a_h} \rix \succ 0$ \\ \hline
Right hyperbola: $\Omega_{Hyp,R}(a_h,b_h)$ & $ \mat{cc} -\frac{R}{a_h} & -Q^{-1} - \frac{J}{b_h}\\ -Q^{-1}+\frac{J}{b_h} & -\frac{R}{a_h} \rix \succ 0$ \\ \hline
Horizontal strip: $\Omega_H(w)$ & $\mat{cc} wQ^{-1} & -J \\ J & wQ^{-1} 
\rix \succ 0$\\ \hline
\end{tabular}
    \caption{DH parametrization of the LMI regions listed in Section~\ref{sec:specificlmis}.}
    \label{tab:specialdhlmi}
\end{table}

We conclude this section with the following remark. 

\begin{remark}[Non-uniqueness of the DH characterization]{\rm
A remark similar to~\cite[Remark 6]{GilS17} also holds for the non-uniqueness of a DH characterization for an $\Omega$-stable matrix. The non-uniqueness of the DH characterization for $\Omega$-stable matrices can be partly characterized by the non-uniqueness of the characteristic function of an LMI region.  For example, a different DH characterization involving two $n \times n$ LMIs is obtained for the horizontal strip in Appendix~\ref{app:s1charcomegah}.
}
\end{remark}

\section{Application of the DH characterization of LMI regions: Nearest $\Omega$-stable matrix problem}\label{sec:appl}

The DH characterization of matrices with eigenvalues in some LMI regions can be used to solve various problems in systems and control. For example,  
it was used to solve the nearest stable matrix problem~\eqref{eq:probdef} for Hurwitz stability in~\cite{GilS17}, 
for Schur stability in~\cite{gillis2019approximating, mamakoukas2020learning}, and for some LMI regions (namely, vertical strips, disks, and conic sectors)  in~\cite{ChoGS20}. It was also used for static-state feedback (SSF) and static-output feedback (SOF) stabilization in~\cite{GilS21}.
Our new general characterizations (see Theorem~\ref{th:genJRQ}) can be used to generalize these algorithms to any LMI region $\Omega$. 
In this section, we illustrate this on the nearest $\Omega$-stable matrix problem.

\subsection{Reformulation}

Let $A\in \R^{n,n}$, $\Omega\subseteq \C$ be an LMI region, and consider the nearest $\Omega$-stable matrix problem~\eqref{eq:probdef}. Then in view of Theorem~\ref{th:genJRQ}, by setting $P=
Q^{-1}$, we can parameterize the set $S_{\Omega}$ of $\Omega$-stable matrices in terms of matrix triplets $(J,R,P^{-1})$ as follows 
\begin{equation*}
    S_\Omega = \left\{(J-R)P^{-1} \in \R^{n,n}:~
    J^T=-J,R^T=R,P \succ 0, \mathcal M_{\Omega}(J,R,P^{-1}) \prec 0
    \right\}.
\end{equation*}  
This leads to the following equivalent formulation of the nearest $\Omega$-stable matrix problem~\eqref{eq:probdef}: 
\begin{equation}\label{eq:reform_stab}
    \inf_{(J,R,P) \in \mathbb DH_{\Omega} } {\|A-(J-R)P^{-1}\|}_F^2, 
\end{equation}
where 
\begin{equation*} 
\mathbb DH_{\Omega}=\left \{ (J,R,P) \in (\R^{n,n})^3:~J^T=-J,\,R^T=R,\,P \succ 0,\,M_{\Omega}(J,R,P^{-1}) \prec 0\right \}.
\end{equation*}
Note that the feasible set in~\eqref{eq:reform_stab} is convex involving only convex LMI constraints, but the objective function is nonconvex.

\subsection{Algorithm}

To solve the nonconvex optimization problem~\eqref{eq:reform_stab}, we follow the same strategy as in~\cite{ChoGS20}. We have adapted the MATLAB code of~\cite{ChoGS20} to handle all the LMI regions summarized in Table~\ref{tab:specialdhlmi}. 
In a few words, the strategy in~\cite{ChoGS20} uses a simple gradient descent method, that is, 
it uses the update 
\[
X \leftarrow \mathcal{P}\big( X - \gamma \nabla F(X) \big), 
\]
where 
$X = (J,R,Q)$ is the variable, 
$F(X) = {\| A - (J-R)P^{-1}\|}_F^2$, 
$\mathcal{P}$ is the projection onto the feasible set $DH_{\Omega}$, which can be computed via semidefinite programming (SDP), 
$\gamma$ is a step size computed using backtracking line search ($\gamma$ is reduced until $F$ decreases), 
and $\nabla F(X)$ is the gradient of $F$ at $X$. 
After each such update, the variables $(J,R)$, with $P$ being fixed, are updated using SDP, which we solve using CVX~\cite{cvx,gb08}. 
An effective and simple initialization for the factors is to choose $P = I_n$ and set $(J,R)$ as the optimal solution of~\eqref{eq:reform_stab} for $P$ fixed.

\subsection{Numerical illustrations}  

As far as we know, only our approach can tackle the nearest $\Omega$-stable matrix problem with LMI regions in full generality. Although Noferini and Poloni~\cite{noferini2020nearest} proposed a general framework for this problem using the Schur decomposition of $A \approx UTU^\top$, they need to project $2 \times 2$ diagonal blocks of $T$ onto the corresponding LMI region. It is unclear how to perform this projection in full generality (in their paper, they work out the details for the cases of Hurwitz and Schur stability). 

To illustrate the use of our algorithm, we consider a similar setting as in~\cite{ChoGS20}: we generate $(J,R,P)$ randomly (namely using \texttt{randn(n)} in
MATLAB), then project it onto the desired LMI region so that $A = (J-R) P^{-1}$ is $\Omega$-stable.  
Then we perturb $A$ using Gaussian noise: 
$A_p = A + N$ 
where $N(i,j) \sim N(0, \sigma)$ where 
$\sigma = \frac{\epsilon{\| A \|}_F }{n}$ where 
$\epsilon > 0$ is a parameter so that the expected value of $\|N\|_F^2$ is equal to $\epsilon^2 {\| A \|}_F^2$ ($\epsilon$ is a measure of the noise-to-signal ratio).  

\paragraph{Example 1} 

\begin{figure}[ht!]
\begin{center}
\includegraphics[width=0.7\textwidth]{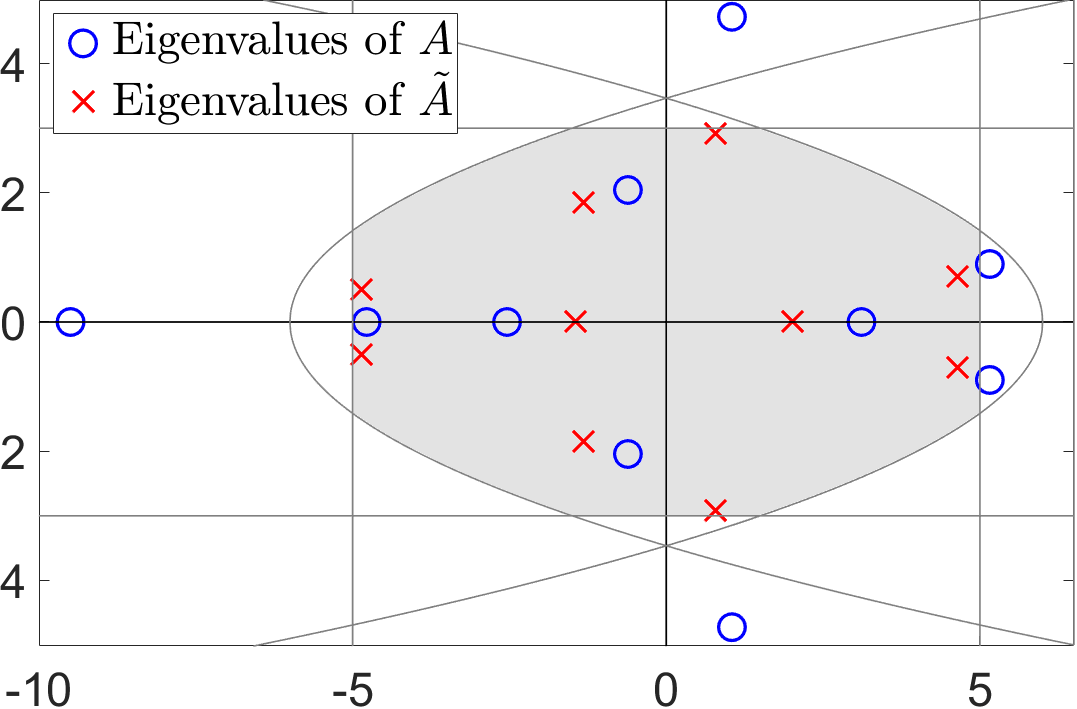}
\caption{Eigenvalues of
$A$, and of its $\Omega$-stable approximation $\tilde{A} = 
(J- R) Q$, where $\Omega$ is the intersection of a vertical strip, a horizontal strip,  and a left and a right parabolic region. \label{fig:ex1}}
\end{center}
\end{figure}

Let us consider $\Omega$ as the intersection of 
\begin{itemize}
    \item a vertical strip between -5 and 5, 
    \item a horizontal strip between -3 and 3,  
    \item  a left parabolic region centered at $(6,0)$ and curvature 1, and 
    \item  a right parabolic region centered at $(-6,0)$ and curvature 1.   
\end{itemize}
We randomly generate a 10-by-10 matrix following the construction described above  with $\epsilon = 1$ (the noise level is high to have several eigenvalues of $A$ outside $\Omega$), and run our algorithm with the identity initialization to approximate $A$ with an $\Omega$-stable matrix $\tilde{A} = (J-R)Q$; see Figure~\ref{fig:ex1} for the illustration. The approximation is computed in approximately 1 minute on a standard laptop, and the relative error is given by $\frac{\|A - \tilde{A} \|}{\|A\|_F} = 18.1\%$.  

\paragraph{Example 2} 

\begin{figure}[ht!] 
\begin{center}
\includegraphics[width=0.7\textwidth]{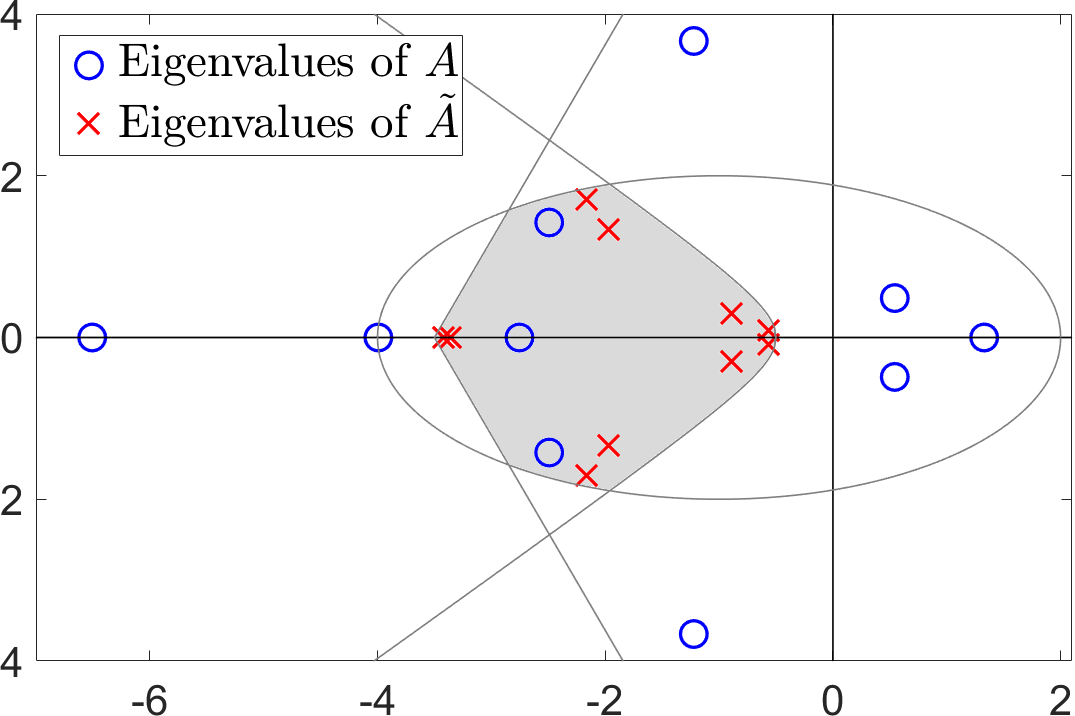}
\caption{Eigenvalues of
$A$, and of its $\Omega$-stable approximation $\tilde{A} = 
(J- R) Q$, where $\Omega$ is the intersection of an ellipse, a right conic sector and a left hyperbolic region. \label{fig:ex2}}
\end{center}
\end{figure}

Let us consider $\Omega$ as the intersection of 
\begin{itemize}
    \item an ellipsoid centered at 
    $(-1,0)$ with horizontal radius 3 and vertical radius of 2, 
 
    \item  a left hyperbolic region centered with semi-majos axis $a_h = b_h = 0.5$, and 
    
    \item  a right conic sector centered at $(-3.5,0)$ with angle $\frac{3}{8} \pi$.   
\end{itemize}
We randomly generate a 10-by-10 matrix following the construction described above with $\epsilon = 1$, and run our algorithm with the identity initialization to approximate $A$ with an $\Omega$-stable matrix $\tilde{A} = (J-R)Q$; see Figure~\ref{fig:ex2} for the illustration. The approximation is computed in approximately 1 minute on a standard laptop, and the relative error is given by $\frac{\|A - \tilde{A} \|}{\|A\|_F} = 24.1\%$.  

\begin{remark}[Code]
The MATLAB code is available  from~\url{https://sites.google.com/site/nicolasgillis/code}. With the code, you can run the examples presented above (we used the random seed 2017 to make the two experiments above reproducible),  
and run any other example.  
\end{remark}

\section{Extended LMI regions for complex matrices} \label{sec:extebdedLMIs}

In this section, we extend the concept of LMI regions to regions that are not symmetric with respect to the real axis but can be represented by LMIs involving complex matrices. We call such regions \emph{extended LMI regions} and define them as follows. 



\begin{definition}[Extended LMI Regions] \label{def:genlmireg} 
A subset $\Omega \subseteq \C$ is called an extended LMI region if there exists 
a Hermitian matrix $B \in \H^{m,m}$ 
and $C \in \C^{m,m}$ such that 
\begin{equation}\label{eq:geneqlmidef1}
\Omega_{\C}(B,C) =\left \{z\in \C:~f_\Omega(z) \prec 0 \right \},
\end{equation}
where 
\begin{equation}\label{eq:geneqlmidef2}
f_\Omega 
\colon 
\mathbb{C} \mapsto \H^{m,m} \; \text{ given by } \; z \mapsto f_\Omega(z) := B + z C + \overline z C^* . 
\end{equation}
\end{definition}

The extended LMI set $\Omega_{\C}(B,C)$ is not in general symmetric with respect to the real line. Further, the rotation, translation and scaling of an LMI region $\Omega$ by nonzero scalar $\alpha \in C$ is in general not an LMI region but an extended LMI region. For example, let $\Omega=\{z \in \C:~B+zC+\overline z C^T \prec 0\}$ be an LMI region as defined in~\eqref{eq:eqlmidef1} and $\alpha \in \C \setminus \{0\}$. If we denote the scaling and rotation of $\Omega$ by $\alpha$ with $\Omega_{\theta}:=\alpha \Omega=\{\alpha z:~z \in \Omega\}$ and the translation of $\Omega$ by $\alpha$ with
$\Omega_{T}:=\Omega+\alpha=\{z+\alpha\in \C: z\in \Omega \}$, then
it is easy to verify that 
$\Omega_{\theta}=\Omega_\C(B,\frac{1}{\alpha}C)$ and $\Omega_T = \Omega_\C(B-\alpha C-\overline \alpha C^*,C)$.

Next, we explain how the main results (Theorems~\ref{Thm:lmicharmain} and~\ref{th:genJRQ}) of Section~\ref{sec:DHcarac} can be generalized to  extended LMI regions by using complex matrices $B$ and $C$ in~\eqref{eq:eqlmidef2}. We first define stability for a complex matrix $A \in \C^{n,n}$.

\begin{definition}
Let $\Omega(B,C) \subseteq \C$ be an extended LMI region~\eqref{eq:geneqlmidef1}. A matrix $A \in \C^{n,n}$ is said to be complex $\Omega$-stable if all eigenvalues of $A$ lie inside $\Omega(B,C)$.
\end{definition}

The following result is a generalization of Theorem~\ref{Thm:lmicharmain} for complex $\Omega$-stable matrices, the proof of which is similar to the proof of Theorem~\ref{Thm:lmicharmain}. We have included the proof for future reference.

\begin{theorem}\label{thm:genlmi}
Let $A \in \C^{n,n}$ and $\Omega(B,C)$ be an extended LMI region in the form~\eqref{eq:geneqlmidef1}. Then $A$ is complex $\Omega$-stable if and only if there exists a Hermitian matrix $X \in \C^{n,n}$ such that 
\begin{eqnarray}\label{eq:gencond}
 X \succ 0 \quad \text{and} \quad \mathcal M_\Omega(A,X)=B \otimes X + C \otimes (AX) + C^* \otimes (AX)^* \prec 0.
\end{eqnarray}
\end{theorem}
\proof ($\Leftarrow$) First suppose that there exists $X \in \C^{n,n}$ satisfying~\eqref{eq:gencond}. Let $\lambda \in \Lambda (A)$ and $x \in \C^{n}\setminus \{0\}$ such that $x^*A=\lambda x^*$. Then 
\begin{eqnarray}\label{eq:proofgenlmi1}
 (I_n \otimes x^*) \mathcal M_\Omega(A,X)  (I_n \otimes x)=f_\Omega(\lambda)x^*Xx.
\end{eqnarray}
This implies that $f_\Omega(\lambda) \prec 0$, since $\mathcal M_\Omega(A,X) \prec 0$ and $X$ is Hermitian positive definite. Thus $\lambda \in \Omega(B,C)$ and hence $A$ is $\Omega$-stable.

($\Rightarrow$) Conversely, let $A$ be $\Omega$-stable. First assume that $A$ is a diagonal matrix (denote it by $\Delta$) with diagonal entries $\lambda_1,\ldots,\lambda_n \in \Omega(B,C)$. Then 
\begin{eqnarray}\label{eq:gencharproof1}
 \mathcal M_\Omega(\Delta,I_n)= \rm{diag}\left (f_\Omega(\lambda_1),\ldots,f_\Omega(\lambda_n)\right ) \prec 0,
\end{eqnarray}
since $A$ is $\Omega$-stable implies that $f_\Omega(\lambda_i) \prec 0$ for all $i=1,\ldots,n$. 
Now suppose $A \in \C^{n,n}$ and let $\Delta$ be the diagonal matrix containing the eigenvalues of $A$ (counting algebraic multiplicities). Let $S$ be an invertible matrix such that  $J=S^{-1}AS$ is the Jordan matrix of $A$. Then there exists a squence of invertible matrices ${T_k}$ such that 
$\lim_{k \rightarrow \infty} T_k^{-1}J T_k = \Delta$. Since $\mathcal M_\Omega(Y,I_n)$ is a continuous function of $Y$, in view of~\eqref{eq:gencharproof1}, we have that 
\[
\lim_{k \rightarrow \infty} \mathcal M_\Omega(T_k^{-1}JT_k,I_n)=\mathcal M_\Omega(\Delta,I_n) \prec 0.
\]
This implies that there exists a positive integer $k$ such that $T:=T_k$ satisfies
$\mathcal M_\Omega(T^{-1}JT,I_n) \prec 0$, or, equivalently $\mathcal M_\Omega(T^{-1}S^{-1}AST,I_n) \prec 0$. Thus, we have 
\begin{eqnarray*}
 (I_n \otimes ST) \mathcal M_\Omega(T^{-1}S^{-1}AST,I_n)  (I_n \otimes {(ST)}^*) = \mathcal M_{\Omega}(A,X) \prec 0,
\end{eqnarray*}
where $X:=(ST)(ST)^*$ is a positive definite matrix, since $S$ and $T$ both are invertible. This shows the existence of $X \succ 0$ such that $M_{\Omega}(A,X) \prec 0$, and hence the proof.
\eproof

The following result is a generalization of Theorem~\ref{th:genJRQ} for the complex matrices that characterizes the set of all complex matrices with eigenvalues in an extended LMI region in terms of complex matrix triplets $(J,R,Q) \in (\C^{n,n})^3$ with Hermitian and definite structures. 

\begin{theorem} \label{th:complexgenJRQ}
Let $\Omega(B,C) \subseteq \C$ be an extended LMI region defined by~\eqref{eq:geneqlmidef1} and let $B,C \in \C^{n,n}$ be the corresponding matrices defined in~\eqref{eq:geneqlmidef2}. Consider a matrix $A \in \C^{n,n}$. Then $A$ is complex $\Omega$-stable if and only if $A=(J-R)Q$
for some $J,R,Q \in \C^{n,n}$ such that $J^*=-J$, $R^*=R$, $Q \succ 0$, and
$\mathcal M_{\Omega}(J,R,Q) \prec 0$, where 
\[
\mathcal M_{\Omega}(J,R,Q) : = B \otimes Q^{-1} +  (C-C^* ) \otimes J - (C + C^* ) \otimes R.
\]
\end{theorem}
\proof The proof follows on the lines of the proof of Theorem~\ref{th:genJRQ} by using Theorem~\ref{thm:genlmi} in place of Theorem~\ref{Thm:lmicharmain}.
\eproof

\small

\bibliographystyle{siam}
\bibliography{ChouGS}

\appendix

\section{DH characterization of the horizontal strip with two 
 $n \times n$ LMIs} 
\label{app:s1charcomegah}


The following theorem provides another DH characterization for the horizontal strip compared to the one provided in Table~\ref{tab:specialdhlmi}.  

\begin{theorem}\label{thm:s1charif}
Let $A=(J-R)Q$ be a DH matrix with $J,R,Q \in \R^{n,n}$ such that $J^T =-J$, $R^
T= R$, and $Q \succ 0$. If 
\begin{equation}\label{eq:s1charhor}
    -wQ^{-1} \prec iJ \prec wQ^{-1},
\end{equation}
then $A$ is $\Omega_H(w)$-stable.
\end{theorem}
\proof
Let $\lambda =\lambda_1+i \lambda_2,\lambda_1,\lambda_2\in \R $ be an eigenvalue of $A$ and let $x \in \C^{n}\setminus \{0\}$ such that $(J-R)Qx=\lambda x$. This implies that $x^*QJQx-x^*QRQx= (\lambda_1+i\lambda_2)x^*Qx$ and thus $ix^*QJQx-ix^*QRQx= (i\lambda_1-\lambda_2)x^*Qx$. By comparing the real and imaginary parts, we have
\begin{equation}\label{eq2:s1charhor}
    \lambda_2=\frac{x^*Q(-iJ)Qx}{x^*Qx},
\end{equation}
since $Q \succ 0$. 
As $J$ satisfies~\eqref{eq:s1charhor}, we have $ -wx^*Qx < x^*Q(iJ)Qx < w x^*Qx$ and thus from~\eqref{eq2:s1charhor} we have that 
\begin{equation}
    -w < \frac{x^*Q(-iJ)Qx}{x^*Qx} < w \quad \implies \quad 
    -w < \lambda_2 < w. 
\end{equation}
This implies $\lambda \in \Omega_H(w)$.
\eproof 

The converse of Theorem~\ref{thm:s1charif} is true when $A$ is  semisimple, that is, algebraic mulitiplicity is equal to the geometric  multiplicity for every eigenvalue $\lambda$ of $A$. More precisely, we have the following result.

\begin{theorem}
Let $A \in \R^{n,n}$ be semisimple $\Omega_H(w)$-stable matrix. Then $A=(J-R)Q$ for some $J,R,Q \in \R^{n,n}$ such that $J^T =-J$, $R^T= R$, $Q \succ 0$, and 
\begin{equation}\label{eq3:s1charhor}
    -wQ^{-1} \prec iJ \prec wQ^{-1}.
\end{equation}
\end{theorem}
\proof In view of Theorem~\ref{th:genJRQ}, A=$(J-R)Q$ for some $J,R,Q \in \R^{n,n}$ such that $J^T =-J$, $R^T= R$, and  $Q \succ 0$, since $A$ is $\Omega_H(w)$-stable. Next, we show that $J$ satisfies~\eqref{eq3:s1charhor}, or equivalently, 
$-wQ \prec Q(iJ)Q \prec wQ$, since $Q\succ 0$. We will prove this by showing that 
\begin{equation}\label{eq4:s1charhor} 
x^*(Q(iJ)Q-wQ)x < 0 \quad \text{and} \quad x^*(Q(iJ)Q+wQ)x>0 \quad \text{for all}~ x \in \C^{n} \setminus \{0\}. 
\end{equation}
As $A$ is semisimple, there exists an orthonormal basis of $\C^{n}$ consisting of eigenvectors of $A$, say $y_1,\ldots,y_n$. Thus, to prove~\eqref{eq4:s1charhor} for every $x$, it is sufficient to show it for vectors $y_j$, $j=1,\ldots,n$. Since $y_j$ is an eigenvector of $A$, we have 
$Ay_j=\lambda_j y_j$ for some eigenvalue $\lambda_j$ of $A$. This implies that 
$(J-R)Qy_j=\lambda_j y_j$ and thus $y_j^*Q(J-R)Qy_j=\lambda_j 
y_j^*Qy_j$. This implies that $\imag (\lambda_j)=\frac{y_j^*Q(-iJ)Qy_j}{y_j^*Qy_j}$ and thus we have 
$-w<\frac{y_j^*Q(-iJ)Qy_j}{y_j^*Qy_j} < w$, since $A$ is $\Omega_H(w)$-stable. 
\eproof 

\end{document}